\documentclass[12pt]{amsart}
\usepackage{epsfig}
\newfont{\sheaf}{eusm10 scaled\magstep1}

\newcommand{\ra}{\ensuremath{\rightarrow}}

\def\eea{\end{eqnarray*}}
\def\bea{\begin{eqnarray*}}

\def\De{{\Delta}}
\def\ga{{\gamma}}
\def\Ga{{\Gamma}}

\newcommand{\Proof}{{\it Proof. }}
\newcommand{\QED}{{\hfill $Q.E.D.$}}

\newtheorem{teo}{Theorem}[section]
\newtheorem{df}[teo]{Definition}
\newtheorem{lem}[teo]{Lemma}

\newtheorem{oss}[teo]{Remark}
\newtheorem{prop}[teo]{Proposition}

\newcommand{\C}{\ensuremath{\mathbb{C}}}

\newcommand{\Q}{\ensuremath{\mathbb{Q}}}
\newcommand{\F}{\ensuremath{\mathbb{F}}}

\newcommand{\N}{\ensuremath{\mathbb{N}}}
\newcommand{\hol}{\ensuremath{\mathcal{O}}}
\newcommand{\HH}{\ensuremath{\mathbb{H}}}
\newcommand{\DD}{\ensuremath{\mathbb{D}}}
\newcommand{\PP}{\ensuremath{\mathbb{P}}}

\newcommand{\I}{\ensuremath{\mathcal{I}}}
\newcommand{\SSS}{\ensuremath{\mathcal{S}}}

\newcommand{\sC}{{\mathcal C}}

\newcommand{\B}{\ensuremath{\mathbb{B}}}

\begin{document}

\title[On varieties whose universal cover is a product of curves  ]
    {On varieties whose universal cover is a product of curves  }

\author{Fabrizio Catanese, Marco Franciosi}


\date{October 10, 2008}

\begin{abstract}       We investigate  a  necessary condition
for a compact  complex
         manifold $X$ of dimension $n$ in order that  its  universal cover be
the Cartesian  product $C^n$ of a curve $C = \PP^1 {\it or} \ \HH$: the
existence of a semispecial tensor $\omega$.

     A semispecial tensor is a  non zero section
        $ 0 \neq \omega \in H^0(X, S^n\Omega^1_X (-K_X) \otimes \eta)  $),
where $\eta$ is an  invertible sheaf of 2-torsion (i.e.,
$\eta^2\cong
\hol_X$).
        We show  that this condition works out nicely, as a sufficient
condition,
when coupled  with some other simple hypothesis, in
    the case of dimension $n= 2$ or $ n= 3$; but it is not sufficient
alone, even in dimension $2$.

In the case of K\"ahler surfaces we use the above results in order
to give a characterization of the surfaces whose universal cover is a
product of
two curves, distinguishing the 6 possible cases.

        \end{abstract}
\maketitle

\section{Introduction}

The beauty of the theory of algebraic curves is deeply related
to the manifold implications of the:

\begin{teo}[Uniformization theorem of Koebe and Poincar\'e]
Let $C$ be a smooth (connected) compact  complex   curve
of genus $g$, and let $  \tilde{C}  $ be its universal cover.
Then
\[
     \tilde{C} \cong
      \left\{ \begin{array}{ll}  {\PP}^1  & \mbox{ if } \  g=0 \\
\C & \mbox{ if } \  g=1 \\
\HH & \mbox{ if } \   g \geq 2 \\
\end{array}
\right.
\]
\end{teo}  ($\HH$ denotes as usual the
Poincar\'e  upper half-plane
$\HH= \{ \tau \in \C :  Im (\tau) > 0\}$, but we shall often refer to it
as the `disk' since it is biholomorphic to
$\DD : = \{ z\in \C : ||z||< 1 \}$).

Hence  a smooth
(connected) compact  complex   curve  $C$  of genus $g \geq 1$
admits a uniformization in the strong sense {\em (ii)} of the following
definition (for $g=0$, only (i) holds):

\begin{df}\label{unif}
A connected complex space $X$ of complex dimension $n$
admits a {\bf Galois
uniformization} if :
\begin{enumerate}
\item[{\em (i)}]
      there is a connected open set $ \Omega \subset \C^n$
and a properly discontinuous group  $ \Ga \subset   Aut (\Omega)$
such that $   \Omega  / \Ga \cong X $
\end{enumerate}

If $X$ is a complex manifold, there is the stronger property
where we require the action of $\Ga$ to be free:
\begin{enumerate}
\item[{\em (ii)}]
      there is a connected open set $ \Omega \subset \C^n$
biholomorphic to the universal cover of $X$
({\bf  strong uniformization}).
\end{enumerate}
\end{df}

Observe that a result of
Fornaess and  Stout  (cf. \cite{f-s}) says that,
if $X$ is an $n$-dimensional complex manifold, then
    there is a connected open set $ \Omega \subset \C^n$
and a surjective holomorphic submersion $ f \colon \Omega \ra X$; i.e.,
every  complex manifold  admits an `\'etale (but not
Galois) uniformization'.

On the contrary,  the condition that the universal cover
be biholomorphic to a bounded domain $ \Omega \subset \subset \C^n$ tends
to be quite exceptional in dimension $ n \geq 2$, where plenty of
simply connected manifolds exist.

    An important remark is that if $ \Omega$ is
bounded  and $\Ga$ acts freely on  $ \Omega$ with compact quotient, then
the  complex manifold $X  : =   \Omega  / \Ga  $  has ample canonical
bundle $K_X$ (see \cite{siegel}): in particular it is a projective
manifold of general type.

Even more exceptional is the case where the universal cover
is biholomorphic to a bounded symmetric domain $ \Omega$,
or where there is a Galois  uniformization
with source a bounded symmetric domain,
and there is already a vast literature on a characterization of these
properties (cf.  \cite{Yau}, \cite{yau1}, \cite{yau}, \cite{Bea}).
The basic result in this direction is S.T.  Yau's uniformization theorem
(explained in \cite{yau1} and \cite{yau}), and for which a very
readable exposition is
contained in the first section of \cite{vz},
emphasizing the role of polystability of the cotangent bundle
for varieties of general type. One would wish nevertheless
for more precise or simple characterizations of the various possible
cases.

The paper  \cite{b-p-t}, which extends work of Yau and Beauville,
especially  \cite{Bea}, gives a nice sufficient condition in order that
the universal cover  of a compact K\"ahler manifold $X$ be biholomorphic
to a product of curves.
If the
tangent bundle
    $T_X$  splits as a sum of line subbundles,
    $T_X= L_1\oplus \cdots \oplus L_n$,  then
      its universal cover $\tilde{X}$
    is  biholomorphic to a product of curves:
    $$\tilde{X} \cong ({\PP}^1  )^r \times  \C^s \times \HH^t ,$$ for  suitable
$r,s,t \in \N$.

The above result is not a characterization, in the sense that the splitting
condition is not a  necessary one, even if we weaken it to the condition that
there is a finite \'etale covering $X' \ra X$ such that the tangent bundle of
$X'$ splits.

The purpose of this  work  is to investigate   to which extent  one can find a
simple characterization of the above property in terms of some necessary and
sufficient    conditions which a compact complex (respectively, K\"ahler)
manifold $X$ must fulfill
     in order that   its universal cover be biholomorphic to a product of
curves.

If we require that the universal cover $\tilde{X} $ be biholomorphic to $
({\PP}^1  )^n$ or $ \HH^n $ we have the following necessary condition
(the case of Kodaira surfaces, cf. \cite{Bea}, shows that $ \tilde{X}
\cong \C^n$
without the K\"ahler assumption does not imply this condition):

\begin{df}
Let $X$ be a complex manifold of complex dimension $n$.

Then a {\bf special tensor} is a non zero section
$ 0 \neq \omega \in H^0(X, S^n\Omega^1_X (-K_X) ) $,
while a {\bf semi special tensor} is a non zero section
$0 \neq \omega \in H^0(X, S^n\Omega^1_X (-K_X) \otimes \eta)$,
     where  $\eta$ is an invertible sheaf
such that $\eta^2\cong \hol_X$.

We shall say that the  semi special tensor  is of unique type
if moreover
$ dim (H^0(X, S^n\Omega^1_X (-K_X) \otimes \eta)) = 1.$
\end{df}

We have in fact:

\begin{prop}\label{nec}
Let
$X$  be a compact complex manifold whose universal
cover is biholomorphic to
$({\PP}^1  )^n $ or to $ \HH^n $: then $X$
admits a semi special tensor.
\end{prop}

As we shall see considering the  two dimensional case,
the existence of a semispecial tensor is   not sufficient in
    order to guarantee a totally split
universal cover, and one has to look for further complementary assumptions,
    one such can be for instance the condition of ampleness of the
canonical divisor $K_X$.

Let us discuss first the case of a  smooth compact complex surface.

Here, a famous uniformization result is the characterization, due to
Miyaoka and Yau, of complex surfaces whose universal cover is the  two
dimensional ball
$\B_2$. It is  given  purely in terms of certain numbers
which are either bimeromorphic or topological invariants.

\begin{teo}[Miyaoka-Yau] Let $X$ be a  compact complex  surface. Then $X
\cong {\B}_2 / \Gamma$ (with $\Gamma$ a cocompact  discrete  subgroup of
$\operatorname{Aut} ({\B}_2)$  acting freely on ${\B}_2$)
     if and  only if
\begin{enumerate}
\item
$K_X^2 = 9 \chi(S) > 0$;
\item  the second plurigenus  $P_2(X) >0$.
\end{enumerate}
\end{teo}

The theorem follows combining Miyaoka's result
(\cite{miy}), that these two conditions imply the ampleness of  $K_X$,
with Yau's
uniformization  result (\cite{Yau}) which proves
the existence of a K\"ahler-Einstein metric.

In the case where $X = ( \HH  \times \HH )/ \Gamma$, with $\Gamma$ a
discrete cocompact subgroup of
     $\operatorname{Aut} (\HH \times \HH) $ acting freely, one has
     $K_X^2 =8 \chi(X)$.

But  Moishezon and Teicher in
     \cite{MT}  showed the existence of  a simply connected
     surface of general  type  (hence with $P_2(X) >0$)
having $K_X^2 =8 \chi(X)$,
so that the above conditions are necessary, but not sufficient.
Our contribution here is a by-product
of our attempt to answer the still open question
     whether  there exists a minimal
     surface of general type with   $p_g(X) =0, K_X^2 =8$
which is not uniformized by
$ \HH  \times \HH$ (one has the same question for $\chi(X) =1, K_X^2
=8$).

    The first result of   this
note is  a  precise characterization of  compact complex  surfaces
      whose  universal cover is the bidisk,
respectively the quadric $\PP^1 \times\PP^1$,
     discussing   whether some
hypotheses can be dispensed with.  We have  the following
result giving a  refinement  of a theorem of S.T. Yau
(theorem 2.5 of \cite{yau}), giving sufficient conditions for  (ii) to hold.

\begin{teo}\label{unibidisk}  Let $X$   be a  compact complex  surface.

   $X$ is strongly uniformized by the bidisk  ( $X \cong  (\HH  \times
\HH)/ \Gamma$ , where $\Gamma$ is a cocompact discrete  subgroup of
$\operatorname{Aut} (\HH \times \HH) $ acting freely ) if and only if

\begin{enumerate}
\item[
{(1*)}] $X$ admits  a semi special tensor of unique type;
\item[
{(2)}] $K_X^2 > 0$;
\item [
{(3)}]  the second plurigenus  $P_2(X)\geq 1$.
        \end{enumerate}

$X$ is biholomorphic to  $ \PP^1 \times \PP^1$ if and only if
\begin{enumerate}
\item[
{(1**)}] $X$ admits  a unique special tensor;
\item[
{(2)}] $K_X^2 = 8$;
\item [
{(3**)}]  the second plurigenus  $P_2(X)= 0$;
\item[
{(4)}] $h^0( \Omega^1_X(- K_X))= 6$
        \end{enumerate}
\end{teo}

In the above theorem one can replace condition (3) by :

$(3^*) P_2(X)\geq 2$,

it is moreover interesting to see that none  of the
above hypotheses can be dispensed with.
The most intriguing examples are  provided by

\begin{prop}\label{elliptic}
There do exist properly elliptic surfaces $X$ satisfying
\begin{itemize}
\item
{ (1)} $X$ admits a   special tensor;
\item
     { (3*)} the second plurigenus  $P_2(X)\geq 2$;
\item
$ q (X) : = dim ( H^1 (\hol_X )) > 0 $;
\item
$K_X^2 = 0$;
\item
$X$ is not birational to a product.
       \end{itemize}
\end{prop}

In this respect, we would like to pose  the following question,
which will be discussed in a later section.

{\bf Question. \ }{\em  Let $X$ be a surface with $ q(X)= 0$
and satisfying {\em (1*)}  and {\em (3*)}: is then $X$ strongly
uniformized by the bidisk?}

Our final result concerning algebraic surfaces whose universal
cover is a product of two curves follows combining the previous Theorem
\ref{unibidisk}  with the following

\begin{teo}\label{surfaces}
Let $S$ be a smooth compact K\"ahler surface $S$.
Then the universal cover of $S$ is biholomorphic to
\begin{enumerate}
\item
$\PP^1 \times \C$ $\Leftrightarrow$ $ P_{12}: = P_{12} (S) = 0$, $q
: = q(S)=1$,
$ K_S^2 = 0$.
\item
$\PP^1 \times \HH$ $\Leftrightarrow$ $P_{12} = 0$,
$q\geq 2$,
$ K_S^2 = 8 (1-q).$
\item
$\C^2$ $\Leftrightarrow$ $P_{12} = 1$, $ q= 1$ or $q=2$, $K_S^2 = 0$.
\item
$\C \times \HH$ $\Leftrightarrow$ $P_{12} \geq 2$, $ e(S) = 0.$

\end{enumerate}
\end{teo}

     Concerning the higher dimensional cases, we restrict our attention
here to the case of
manifolds with ample canonical divisor $K_X$ which, by Yau' s theorem
(\cite{Yau}) admit  a canonical K\"ahler-Einstein metric.

    Assume now that $X$ admits a semi special tensor
    $\omega
\in H^0(X, S^n\Omega^1_X (-K_X) \otimes \eta)$.  Then
by \cite[p.272]{yau1} and  by \cite[p.479]{yau}   (see also \cite[p.10]{vz})
$\omega$ induces on the   tangent bundle $T_X$ a homogeneous
    hypersurface  $F_X$ of relative degree $n$ which is parallel with
respect to  the
K\"ahler-Einstein metric.

    In particular,  take a point $x\in X$,  and consider the hypersurface
    of the projectivized tangent bundle induced by $F_X$:
its  fibre over $x$  is a   projective
    hypersurface $F_{X,x}$ of degree $n$
     which is invariant for the action of the (restricted) holonomy group
    $H \subset U(n)$ ($H$ is the connected component of the identity in
the holonomy group).

In this situation, assume that we can prove (possibly passing to a
finite \'etale
covering of $X$) that the holonomy leaves invariant a complete flag.
Then, since the holonomy is unitary, it follows that $H \subset U(1)^n$
and we can conclude, either by Berger's  classical theorem (\cite{berger}),
or by \cite{b-p-t}, that   the universal cover of $X$ turns out to be
$\HH^n$.

     In the three dimensional case the existence of a special tensor is
enough in order to guarantee
     such a splitting.

     \begin{teo}\label{3dim}  Let $X$ be a compact complex manifold
of dimension
$ n \leq 3$.
    Then the following two conditions:
\begin{enumerate}
\item[{ (1)}]
$X$ admits a semi special tensor;
\item[{ (2*)}]
$K_X $ is ample
        \end{enumerate}
hold  if and only
if   $X \cong  (\HH^n  )/ \Gamma$   (where $\Gamma$
is a cocompact discrete  subgroup of $\operatorname{Aut} (\HH^n ) $ acting
freely ).
\end{teo}

In dimension $\geq 4$, the above conditions are no longer sufficient.
The natural category  which is relevant   to consider  is the category of
Hermitian symmetric spaces  of noncompact
type,  since  by  the
    theorem of Berger-Simons  an
irreducible (in the sense of  De Rham's  theorem)
K\"ahler   manifold $X$  of dimension $n$ with ample canonical 
divisor $K_X$  has
holonomy
$H
\neq U(n)$ if and only if $X$ is a Hermitian symmetric space
     of rank $\geq 2$ (see \cite{yau1}, and \cite{vz}, section 1, page 300).

One has   the Cartan realization of a
Hermitian symmetric space of noncompact type as a
bounded symmetric domain, and    by
the classical result of Borel on  compact Clifford-Klein forms
    (see \cite{bo63})  any
    bounded symmetric domain $X$ of dimension $n$
    admits a compact complex analytic Clifford-Klein
     form,  that is  a compact complex manifold $X'$
whose universal  covering is isomorphic to
     $X$.

The above  results translate the question whether a compact
complex manifold $X$ admitting a semi special tensor and with
ample canonical divisor $K_X$ has the polydisk as universal cover
into a purely Lie theoretic problem, the problem of existence
of holonomy invariant hypersurfaces of degree $n$.

We leave aside for the moment this more general investigation,
for which some partial results are contained in the appendix,
due to A.J. Di Scala, who answered some of our questions.

     For the bounded domain  $ \Omega \subset \C^4 \cong Mat (2,2,\C) : =
M_{2,2}(\C)$,
     $ \Omega =     \{ Z \in  M_{2,2}(\C) \  :  \  \operatorname{I_2} - ^{t} Z
    \cdot \overline{Z} >0\}$, the Cartan  realization of the Hermitian
symmetric space
     $SU(2,2) / S (U(2) \times U(2))$, Di Scala  pointed out that the holonomy
action of $  (A, D) \in S (U(2) \times U(2)$ is given
by $ Z \mapsto A Z D^{-1}$. Hence the square of the
determinant yields an invariant hypersurface of degree $4$ which
is twice a
smooth quadric (and this is indeed the only other possible case).

Using this simple but important observation, we get the following

\begin{teo}\label{4dim} There exist  compact K\"ahler manifolds $X$,   for
each dimension $ n \geq 4$, such that
\begin{enumerate}
\item
[{ (1)}] $X$ admits a   special tensor;
\item
     [{ (2*)}]
$K_X $ is ample
\end{enumerate}
and whose universal cover  $\tilde{X}$  is not $ \cong \HH^n$ (i.e., is not
a product of curves).
\end{teo}

\section{Preliminaries and remarks}\label{preliminari}

\subsection{ Notation.}
      $X$  denotes throughout the paper a  smooth compact complex manifold
of dimension $n$.

We use  the standard
notation of algebraic geometry: $\Omega^1_X$  is the cotangent
bundle (locally free sheaf),
      $T_X$ is the holomorphic tangent bundle,
$c_1(X)$, $c_2(X)$ are the Chern classes of $X$.
      $K_X$ is a canonical divisor on $X$, i.e., $\Omega^n_X= \hol_X(K_X)$ and
the  $m$-th plurigenus is defined as $P_m (X) :=
h^0(X, m K_X)$.

In particular, for $m=1$, we have
       the geometric genus of $X$ $p_g(X): = h^0(X, K_X)$,
while $q (X): = h^1(X, \hol_X)$ is classically called
the {\bf irregularity } of $X$.

Finally,
      $\chi (X) : = \chi (\hol_X) $
is the holomorphic Euler Poincar\'e characteristic of $X$,
whereas  $e(X)$ denotes the topological Euler Poincar\'e
characteristic of $X$.

In the surface case ($n=2$), $\chi (X)= 1 +p_g (X) - q (X)$.

      With a slight abuse of notation,  we do not distinguish between
invertible sheaves, line
bundles and divisors, while the symbol
     $\equiv$ denotes  linear equivalence of divisors.

\subsection{ Necessary conditions.}

\hfill\break
First of all
notice that the existence of a
semi special tensor corresponds to the existence of a special tensor
on an \'etale  double
cover of our manifold:
\begin{oss}\label{double}
A complex manifold $X$  admits a semi special tensor if and only if
it has an unramified  cover $X'$ of degree at most two which admits a special
tensor.
\end{oss}
\begin{proof}

Assume that we have an
invertible sheaf $\eta$ such that $\eta^2\cong \hol_X$,
     $\eta \not\cong \hol_X$. Take the corresponding double connected
      \'{e}tale covering  $ \pi :  X' \ra X$ such that
$\pi _* \hol_{  X' }  \cong \hol_{  X } \oplus \eta  $ and observe that
      $$  H^0(X', S^n\Omega^1_{X'} (-K_{X'})) \cong
H^0(X, S^n\Omega^1_{X} (-K_{X})) \oplus H^0(X, S^n\Omega^1_{X}
(-K_{X}) \otimes \eta).$$
Whence, there is a  special tensor on $X'$ if and only if there is
a semi special tensor on $X$.
       \end{proof}

Let us now show that  if  $X$ is isomorphic to $ ({\PP}^1  )^m / \Ga$ or $
(\HH )^m / \Ga$  then $X$ admits  a semi special tensor.

\hfill\break
{\em Proof of Prop. \ref{nec} } .
    Let us remark first that for  a simply connected curve $C$ , with
    $C \cong  \PP_1$ , or $C \cong   \HH$,
    and any integer $m$,
    the group of automorphism of $C^m$,
$\operatorname{Aut} (C^m) $, is the semidirect product
     of $(\operatorname{Aut} (C))^m$ with the symmetric group
${\mathfrak S}_m$, hence for every subgroup $\Ga_C$
of  $\operatorname{Aut} (C^n) $ we have a diagram:
        \[
      \begin{array}{ccccccc}
       1 \ra &  (\operatorname{Aut} (C))^m  & \ra & \operatorname{Aut}
(C^m) &  \ra {\mathfrak S}_m &  \ra 1\\
       &     \bigcup & & \bigcup &   \bigcup & \\
      1 \ra & \Gamma^{0}_C & \hookrightarrow & \Gamma_C  & \ra   H_C & \ra 1
.\\
       \end{array}
        \]
Let now $X \cong (C^n) / \Ga$ be a compact complex manifold whose
universal cover $\tilde{X}$ is isomorphic to  $C^n $.
Then $X$ admits a
semi special tensor, induced by the following special tensor:

$$ \tilde{\omega} : =
\frac{ \operatorname{d}z_1 \otimes \dots \otimes
\operatorname{d}z_n}{\operatorname{d}z_1\wedge \dots \wedge
     \operatorname{d}z_n},
     $$
where $(z_1, \dots, z_n)$ is the standard system of coordinates on $C= \HH^n
$, respectively on the standard open set
     $\C^n \subset ({\PP}^1)^n$ (observe that $\tilde{\omega}$ is in this case
everywhere regular).

$\tilde{\omega}$ is clearly
invariant for $(\operatorname{Aut} (C))^n$ and for the alternating
subgroup ${\mathfrak A}_n.$ Let $\eta$ be the    2-torsion invertible
sheaf   on
$X$ associated to the signature character of  ${\mathfrak S}_n$ restricted to
$H_C$:
then clearly $\tilde{\omega}$ descends to a semi special tensor $\omega
\in H^0(X, S^n\Omega^1_X (-K_X) \otimes \eta)$.

\hfill$\Box$

In the more general case where the universal cover is a product of curves,
we have the following proposition:

\begin{prop}\label{fibrebundle}
We have a
homomorphism
$$\Phi: \operatorname{Aut}(({\PP}^1  )^r \times  \C^s \times \HH^t)
\ra  \operatorname{Aut}( \C^s \times \HH^t)$$
which is injective on any subgroup $\Ga$ which acts freely. Moreover,
   if
$\Ga_2\subset \operatorname{Aut}( \C^s \times \HH^t) $ is the image of
$\Ga$ under $\Phi$, $\Ga_2$ acts also freely, and $\Ga_2$ acts properly
discontinuosly if $\Ga$ is properly discontinuos.

In particular, if
    $X \cong (({\PP}^1  )^r \times  \C^s \times \HH^t ) / \Ga$, with $\Ga$ a
cocompact discrete  subgroup of $\operatorname{Aut}(({\PP}^1  )^r \times
\C^s \times \HH^t)$ which acts freely, then
    the natural projection $$(({\PP}^1  )^r \times  \C^s \times \HH^t) / \Ga
\ra    (\C^s \times \HH^t) / \Ga_2 $$
     inherits a
$({\PP}^1  )^r-$bundle structure.

\end{prop}
Before giving the proof let us point out the following:

\begin{lem}\label{P1} Let $\psi \in\operatorname{Aut}(({\PP}^1  )^r)$
be an
automorphism. Then $\psi$ has a fixed point.

\end{lem}
\begin{proof} For $r=1$ this is well known, since there exists an
eigenvector for each $ A \in GL(2,\C)$.

For $ r \geq 2$ any   automorphism $\psi
\in\operatorname{Aut}(({\PP}^1  )^r)$ is of the form
$$ (\psi (x))_i =  \psi_i (x_{\sigma(i)})$$
for a suitable permutation $\sigma$ of $\{1, \dots r\}$.
Therefore a fixed point is a solution to the system of equations
$$ x_i =  \psi_i (x_{\sigma(i)})  \  (i = 1, \dots r ).$$
Using the cycle decomposition of  $\sigma$ we easily reduce to the case
where $\sigma = (1,2,\dots r)$ and it suffices to find   a
solution to $ x_1 =  \psi_1 \circ \dots \psi_r (x_1)   $.

\end{proof}

\begin{proof} {\em of Prop \ref{fibrebundle}}.
Let
    $\phi\in \operatorname{Aut}(({\PP}^1  )^r \times  \C^s \times \HH^t)$.

Let $\Phi_2$ be the composition $p_2 \circ \phi$, where

$p_2 : ({\PP}^1  )^r \times  \C^s \times \HH^t \ra
\C^s \times \HH^t$

is the second projection.

Now,  for every point $p \in
\C^s
\times
\HH^t$,
$\Phi_2$ is constant on $({\PP}^1  )^r \times \{p\}$ since $({\PP}^1  )^r$
is compact.

    Hence  $\phi$ induces $\phi_2 \in \operatorname{Aut}(
     \C^s \times \HH^t)$.

Assume that $\phi $ acts freely, and that $\phi_2 $ has a fixed point $p$.
Then
$\phi$ acts on
$({\PP}^1  )^r \times \{p\}$ and it has a fixed point there by the previous
lemma: whence
$\phi$ is the identity.

If the action of $\Ga$ is properly discontinuous, then for any compact $K
\subset (\C^s \times \HH^t)$, also $({\PP}^1  )^r \times  K$ is compact;
hence the set $ \Ga_2(K,K) = \Ga(({\PP}^1  )^r \times K,({\PP}^1  )^r
\times K)$ is finite. Therefore $\Ga_2$ is also properly discontinuous.

\end{proof}
\begin{oss}\label{Inoue}
We also have a
homomorphism
$$\Phi: \operatorname{Aut}(({\C}^1  )^r \times   \HH^t)
\ra  \operatorname{Aut}(  \HH^t)$$
However, as shown by the case of Inoue surfaces, if
    $X \cong ((  \C^r \times \HH^t ) / \Ga$, where
    $\Ga$ is a
cocompact discrete  subgroup of $\operatorname{Aut}(({\PP}^1  )^r \times
\C^s \times \HH^t)$ which acts freely,
then the image group
$\Ga_2\subset \operatorname{Aut}(  \HH^t) $ does not necessarily
act properly discontinuously. One needs for this the
assumption that $X$ be K\"ahler.
\end{oss}


\section{Surfaces whose universal cover is  a product of curves}

In the case of surfaces the existence of a special tensor, as we are now
going to explain, is equivalent to the existence of a trace zero
endomorphism of the tangent bundle: and  if this endomorphism
is not nilpotent, one obtains a splitting of the tangent bundle.

Let us recall a result of  Beauville which
characterizes  compact complex surfaces whose  universal
cover is a product  of
two complex curves (cf. \cite[Thm. C]{Bea}).

\begin{teo}[Beauville] \label{teo:beauville}
    Let $X $ be a
compact  complex surface.
    The tangent bundle
$T_X$ splits as a
direct sum of two line bundles
     if and only if either $X$ is a
special Hopf surface
or the universal covering space
of
$X$  is a
product
$U\times V$  of two complex curves and the group $\pi_1(X)$
acts
diagonally on
$U \times V$.
\end{teo}

Given a direct
sum
decomposition of the cotangent bundle $\Omega^1_X \cong L_1\oplus
L_2$,
    Beauville shows   moreover that $(L_1)^2 =( L_2)^2=0$   (cf.
\cite[4.1,
4.2]{Bea}) hence
    \[
    K_X \equiv L_1 +L_2 \hspace{ 1 cm}
c_1(X)^2  = 2\cdot (L_1 \cdot L_2 )
= 2\cdot c_2(X), \  i.e.,  K_X^2 = 8 \chi(X).
    \]

Let us now consider
the  bundle
$\operatorname{End} (T_X)$
    of endomorphisms of the
tangent bundle.   We can write
$\operatorname{End} (T_X) = \Omega^1_X
\otimes T_X $  and  since from the
nondegenerate bilinear  map
$$
\Omega^1_X \times \Omega^1_X \longrightarrow \Omega^2_X \cong
K_X$$
we get $T_X=  (\Omega^1_X )^{\vee} \cong \Omega^1_X (-K_X)$
   we  have an isomorphism

$\operatorname{End}
(T_X)\cong  \Omega^1_X \otimes \Omega^1_X
(-K_X).$

   Let us see how
this isomorphism
works  in local
coordinates $(z_1, z_2)$. I.e.,
   let
us see how an element
$\frac{ \operatorname{d}z_i
\otimes
\operatorname{d}z_j}{\operatorname{d}z_1 \wedge
\operatorname{d}z_2}$
   in $ \Omega^1_X \otimes \Omega^1_X (-K_X)$
acts on a vector of  the
form
   $\frac{\partial}{\partial z_h}$.  We
have
   $$\frac{ \operatorname{d}z_i
\otimes
\operatorname{d}z_j}{\operatorname{d}z_1
\wedge
\operatorname{d}z_2}  \bigl( \frac{\partial}{\partial z_h}
\bigr) =
\left\{ \begin{array}{cl}
   \frac{
\operatorname{d}z_j}{\operatorname{d}z_1 \wedge
\operatorname{d}z_2} & \mbox{ if } h=i  \\
0 & \mbox{ if } h\neq i
\\
\end{array} \right.
$$
   In turn,
$\displaystyle{\frac{\operatorname{d}z_j
}{\operatorname{d}z_1
\wedge
   \operatorname{d}z_2} }$ evaluated on
$\operatorname{d}z_k$ gives
$\displaystyle{\frac{\operatorname{d}z_j
\wedge
    \operatorname{d}z_k}{\operatorname{d}z_1 \wedge
\operatorname{d}z_2} }$.

Therefore  a generic element  $
\displaystyle{ \sum_{i,j} a_{ij}
\frac{ \operatorname{d}z_i
\otimes
\operatorname{d}z_j}{\operatorname{d}z_1 \wedge
\operatorname{d}z_2}}$
    corresponds to an endomorphism, which,  with
respect to the  basis
$\bigl\{ \frac{\partial}{\partial
z_1},\frac{\partial}{\partial z_2}
\bigr\}$ is expressed by the
matrix
$$ \begin{pmatrix} -a_{12} & -a_{22} \\ a_{11} & a_{21}
\\
\end{pmatrix}
$$  In particular for  the symmetric  tensors (i.e.,
$a_{12}= a_{21}$),
respectively for the skewsymmetric tensors  (i.e.,
$a_{12}=-a_{21},
a_{11}=a_{22}=0$) the following isomorphisms
hold:
$$ S^2( \Omega^1_X ) (-K_X) \cong \bigg\{
\begin{pmatrix} -a &
-a_{22} \\ a_{11} & a \\
\end{pmatrix}   \bigg\} \ ; \ \  \hspace{
0,5 cm}
    {\bigwedge}^2( \Omega^1_X ) (-K_X) \cong
\bigg\{
\begin{pmatrix} b & 0 \\ 0 & b \\
\end{pmatrix}
\bigg\}
$$
We can summarize the above discussion in the
following

\begin{lem}\label{split}
If $X$ is a complex surface there
is a natural isomorphism
between the sheaf $S^2( \Omega^1_X )
(-K_X)$
and the sheaf of trace zero endomorphisms of the (co)tangent
sheaf
$\operatorname{End}^0 (T_X)
\cong \operatorname{End}^0
(\Omega^1_X)$.

A special tensor $\omega \in H^0 (S^2( \Omega^1_X ) (-K_X))$
with
nonzero determinant $ det (\omega) \in \C$ yields an
eigenbundle
splitting
$\Omega^1_X\cong L_1 \bigoplus L_2$ of the
cotangent
bundle.

If instead $ det (\omega) = 0 \in \C$,  the
corresponding endomorphism
$\epsilon$ is nilpotent and yields an
exact sequence of sheaves
    \[
    0 \ra L \ra \Omega^1_X \ra
{\mathcal{ I}}_Z L (-\Delta) \ra 0
     \]
where $L : = ker (
\epsilon)$ is invertible, $\Delta$ is an effective divisor,
and $Z$
is a 0-dimensional subscheme(which is a local complete
intersection).

We have in particular
$K_X \equiv 2L - \Delta$
and
$c_2(X)= {length}(Z) +
    L\cdot
(L-\Delta)$.

\end{lem}

\Proof
We need only to observe that $ det
(\omega)$ is a constant,
since  $det (\operatorname{End} (T_X)) =
det (\operatorname{End}
(\Omega^1_X)) \cong \hol_X$.

If $ det
(\omega) \neq 0$, there is a constant $c \in \C \setminus \{
0\}$
such that  $ det (\omega) = c^2$, hence at every point of $X$
the
endomorphism $\epsilon$ corresponding to the special tensor $
\omega$
has two distinct eigenvalues $ \pm c$.

Let $\omega \in
H^0(S^2\Omega^1_X (-K_X) )$, $\omega \neq 0$,  be
   such that $\det
(\omega)=0$.
   Then the corresponding endomorphism $ \epsilon$
   is
nilpotent of order 2, and there exists an open nonempty subset
$U
\subseteq X$ such that
   $\operatorname{Ker}(\epsilon_{|U})  =
\operatorname{Im}(\epsilon_{|U})$.
   At a point $p$ where
$\operatorname{rank} (\epsilon)=0$, in local
coordinates the
endomorphism $\epsilon$ may be expressed by
   $$ \begin{pmatrix} a &
b \\ c & -a \\
\end{pmatrix}  \ \ a,b,c \mbox{ regular functions
such that } a^2 = - b\cdot c
$$
   Let $\delta :=
\operatorname{G.C.D.} (a,b,c)$. After dividing by
$\delta$, every
prime factor
   of $a$ is either not in $b$, or not in $c$, thus we
can write
   $$ -b= \beta^2 \hspace{ 1 cm}  c=\gamma^2 \hspace{ 1 cm}
a = \beta\cdot
\gamma $$
   Therefore  we obtain
   $$  \begin{pmatrix}
u\\ v \\
\end{pmatrix}  \in \operatorname{Ker}{\epsilon}
\Longleftrightarrow
\left\{\begin{array}{l}  a\cdot u + b\cdot v=0
\\ c\cdot u - a\cdot v =0
\end{array} \right. \Longleftrightarrow
\gamma \cdot u - \beta \cdot v =0
\Longleftrightarrow
\begin{pmatrix} u\\ v \\
\end{pmatrix} = \begin{pmatrix}
\beta
\cdot f\\
\gamma\cdot f \\
\end{pmatrix}
$$
    and, writing our
endomorphism $\epsilon$ as $\epsilon =  \delta \cdot
\alpha$, we
have
   $$ \operatorname{Im}(\alpha) =
   \left\{\begin{array}{l}
\beta
\cdot \gamma \cdot u - \beta^2 \cdot v = \beta \cdot(\gamma
\cdot u -
\beta \cdot v)\\
\gamma^2 \cdot u - \gamma \cdot \beta
\cdot v =
\gamma \cdot (\gamma \cdot
u - \beta \cdot v)
\end{array}
\right.
   $$

   Let $Z$ be the 0-dimensional scheme defined by $\{
\beta=\gamma=0\}$ and
$\Delta$ be the Cartier divisor defined
   by
$\{\delta=0\}$.

   From the above description  we deduce that the
kernel of $\epsilon$
is a line bundle $L$ which fits in the
following
exact sequence:
   \[
   0 \ra L \ra \Omega^1_X \ra
{\mathcal{ I}}_Z L
(-\Delta) \ra 0.
    \]
   Taking the total Chern
classes we infer that:
$K_X \equiv 2L - \Delta$ as divisors
on $X$
and
$c_2(X)= {length}(Z)
+
   L\cdot (L-\Delta)$.
\qed

\begin{lem}\label{blowup}
Let $X$ be a
complex surface and
let $X'$ be the blow up of $X$
at a point
$p$.
Then a special tensor
$\omega'$ on $X'$ induces a special tensor
$\omega$ on $X$,
and the
converse only holds if and only if $\omega$
vanishes at $p$
(in
particular, it must hold : $det (\omega) =
0$).
\end{lem}

\Proof
First of all, $\omega'$  induces a special
tensor  on $ X \setminus \{p\}$, and
by Hartogs'  theorem the latter
extends to a   special tensor
$\omega$ on $X$.

Conversely, choose
local coordinates $ (x,y)$ for $X$ around $p$
and take a local
chart
of the blow up with coordinates $ (x,u)$
where $y=u x$.

Locally
around $p$ we can write
$$
\omega =
\frac{a(\operatorname{d}x)^2 + b
(\operatorname{d}y)^2+
c
(\operatorname{d}x\operatorname{d}y)}{\operatorname{d}x\wedge
\operatorname{d}y}
$$
The pull back $\omega'$ of $\omega$ is given by the following
expression:
$$
\frac{a(\operatorname{d}x)^2 + b (u
\operatorname{d}x  +
x \operatorname{d}u)^2 + c
(u\operatorname{d}x+x\operatorname{d}u)
\operatorname{d}x}
{x\operatorname{d}x\wedge
\operatorname{d}u}
=$$
$$
= \frac{\operatorname{d}x^2
( a+bu^2+cu) +
bx^2\operatorname{d}u^2 +
( 2 bu x + c x )
\operatorname{d}x\operatorname{d}u}
{x\operatorname{d}x\wedge
\operatorname{d}u}
,$$
hence $\omega'$ is regular if and only if
$ \frac{ a+bu^2+cu}
{x}$ is a regular function.

This is obvious if $a,b,c$ vanish at
$p$, since then their pull back
is divisible by $x$. Assume on the
other side that $a,b,c$ are constant:
then we get a rational
function
which is only regular if
$a = b = c =
0.$
\qed

\begin{lem}\label{rational}
Let $X$ be a compact minimal
rational surface admitting a special tensor
$\omega$. Then $ X \cong
\PP^1 \times \PP^1$ or $ X \cong
   \F _n , n \geq 2$. If moreover the special tensor is unique,
then   $ X \cong
\PP^1 \times \PP^1$ or $ X \cong
   \F _2$.
\end{lem}

\Proof
Assume that $X$ is a $\PP^1$ bundle over a
curve $B \cong
\PP^1$, i.e., a ruled surface $ \F _n$ with $ n \geq
0$. Let $\pi \colon X \ra B$
the projection.

By the exact sequence

$$ 0 \ra \pi^* \Omega^1_B \ra \Omega^1_X \ra \Omega^1_{X|B} \ra
0$$
and since on a general fibre $F$ the subsheaf $\pi^* \Omega^1_B$
is trivial,
while the quotient sheaf $\Omega^1_{X|B}$ is negative, we
conclude that any
endomorphism $\epsilon$ carries $\pi^* \Omega^1_B$
to itself.
If it has non zero determinant we can conclude by Theorem
\ref{teo:beauville}
that $ X \cong \PP^1 \times \PP^1$. Otherwise,
$\epsilon$ is nilpotent
and we have a nonzero element in $ {\rm
Hom}(\Omega^1_{X|B}, \pi^*
\Omega^1_B)$.

Since these are invertible
sheaves, it suffices to see when
$$ H^0 ( \hol_X ( 2 \pi ^* K_B - K_X
)) \neq 0.  $$
But, letting $\Sigma$ be the section with
selfintersection
$\Sigma^2 = -n$, our vector space equals
$ H^0 (
\hol_X ( 2 \Sigma + (n -2) F )) .  $
Intersecting this divisor  with
$\Sigma$ we see that (since each time the intersection number
with
$\Sigma$ is negative)  $ H^0 ( \hol_X ( 2 \Sigma + (n -2) F )) = H^0
(\hol_X (\Sigma + (n -2) F ))=  H^0 ( \hol_X ( + (n -2) F )) $.
This space has dimension $n-1$, whence our claim follows for the
surfaces $\F_n$.

There remains the case where $X$ is $\PP^2$.

In this case
$\epsilon$ must be a nilpotent endomorphism by
Theorem
\ref{teo:beauville}, and it cannot vanish at any point by our
previous result
on $ \F _1$. Therefore the rank of  $\epsilon$
equals 1 at each point.
By lemma \ref{split} it follows that there is
a divisor $L$ such that $K_X = 2 L$,
a
contradiction.
\qed

\subsection{Proof of Theorem \ref{unibidisk}}

\begin{proof}
If $X$ is strongly uniformized by the
bidisk, then $K_X$ is ample,
in particular $K^2_X \geq 1 $ and, since
by Castelnuovo's theorem
$ \chi(X) \geq 1$,
by the vanishing theorem
of Kodaira and
Mumford it follows that
$ P_2(X) \geq 2$ (see
\cite{cm}).

Thus one
direction follows from proposition  \ref{nec},
except that we shall
show only later that (1*) holds.

Assume conversely that $(1), (2)$
hold.
Without loss of generality we may assume by lemma \ref{blowup}
that $X$ is
minimal, since $K_X^2$ can only decrease via a blowup and
the bigenus is a
birational invariant.

$K^2_X \geq 1 $ implies that
either the
surface
$X$ is of general type, or it is a rational
surface.

These two cases are  distinguished by the
respective properties  (3) (obviously  implied by (3*)),
guaranteeing that $X$ is of general type,
and (3**) ensuring that $X$ is rational.

Let us first  assume  that  $X$ is of general type
and, passing to
an \'etale double cover if necessary, that
$X$ admits a special
tensor.

By the cited Theorem \ref{teo:beauville} of \cite{Bea} it
suffices
to find a decomposition of the cotangent bundle
$\Omega^1_X$
as a direct sum of two  line bundles $L_1 $ and $L_2$.

The two line
bundles $L_1$, $L_2$ will be given as  eigenbundles of a
diagonizable
endomorphism $\epsilon \in \operatorname{End} (\Omega^1_X)$.

Our
previous discussion shows then that it is sufficient to show that
any
special tensor cannot yield a nilpotent endomorphism.

Otherwise,
by lemma \ref{split}, we can write
    $2L \equiv K_X + \Delta$ and
then deduce that $L$ is a big divisor since
$\Delta$ is
    effective
by construction and $K_X$ is big because  $X$ is of general
type.
This assertion gives the required contradiction since  by
    the
Bogomolov-Castelnuovo-de Franchis Theorem (cf. \cite{Bog})  for
an invertible subsheaf $L$  of  $\Omega^1_X$ it is
$h^0(X, mL) \leq
O(m) $, contradicting the bigness of $L$.

There remains to show
(1*).  But if   $h^0(X, S^2\Omega^1_X (-K_X) ) \geq 2$
then, given a
point $ p \in X$, there is a special tensor which is
not invertible
in
$p$, hence a special tensor with vanishing determinant, a
contradiction.

   If $X$ is a rational surface we use the hypothesis $K_X^2 = 8$,
ensuring that $X$ is a surface $\F_n$; then, by
lemma
\ref{rational} we conclude that either $ X \cong
\PP^1 \times \PP^1$ or $ X \cong
   \F _2$. In the former case $ h^0 (\Omega^1_X (-K_X)) = 6$,
in the latter case $ h^0 (\Omega^1_X (-K_X)) = 7$.
\end{proof}

\section{Elliptic surfaces
with a
special tensor
   not birational to a product of curves}


\hfill\break
In
this section we are going to prove proposition \ref{elliptic}.

We
consider surfaces $X$ with bigenus $P_2(X) \geq 2$
(property (3*)),
therefore their Kodaira dimension equals 1 or 2,
hence either they
are properly (canonically) elliptic, or they are of
general
type.

Since we took already care of the latter case in the main
theorem
\ref{unibidisk},
we restrict our attention
here to the former case, and try to see when does
a properly elliptic
surface admit a special tensor (we can reduce to
this situation
in
view of remark \ref{double}). We can moreover assume that the
associated
endomorphism $\epsilon$ is nilpotent by theorem
\ref{teo:beauville}.

Again without loss of generality we may assume
that $X$ is minimal by
virtue of lemma \ref{blowup}.

	\Proof
Let
$X$ be a minimal properly elliptic surface and let $f : X \ra B$
be
its (multi)canonical elliptic fibration.
Write any fibre $f^{-1} (p)
$ as $F_p = \sum_{i=1}^{h_p} m_i  C_i$
and, setting $ n_p : = G.C.D.
(m_i)$, $ F_p = n_p F'_p$,
we say that a fibre is multiple if $n_p >
1$.
By Kodaira's classification (\cite{Kodaira}) of the singular
fibres
we know that in this case $ m_i = n_p ,  \forall i.$

Assume
that the multiple fibres of the elliptic fibration are $n_1 F_1',
\dots
, n_r F_r'$, and consider the divisorial part of the critical
locus
$$ \SSS_p : = \sum_{i=1}^{h_p} (m_i-1 )  C_i , \ \  \SSS : =
\sum_{p
\in B}\SSS_p $$
    so that we have then the exact sequence
$$
0 \ra f^* \Omega^1_B (  \SSS ) \ra
\Omega^1_X \ra  \I_{\sC}\
\omega_{X|B} \ra 0,$$
where $\sC$ is a 0-dimensional (l.c.i.)
subscheme.

For further calculations we separate the divisorial part
of the critical locus
as the sum of two disjoint effective
divisors,
the multiple fibre contribution and the rest:
$$ \SSS_m  :
= \sum_{i=1}^r
(n_i -1 ) F'_i, \  \hat {\SSS} : = \SSS  - \SSS_m
.$$

Let us assume that we have a nilpotent endomorphism
corresponding to another
exact sequence
\[
    0 \ra L \ra \Omega^1_X
\ra {\mathcal{ I}}_Z L (-\Delta) \ra 0,
     \]
in turn determined by a
homomorphism
$$ \epsilon' :  {\mathcal{ I}}_Z L (-\Delta) \ra L,$$
i.e., by a section $$s \in H^0 (\hol_X (\Delta)) = $$
$$=H^0 (\hol_X (2 L - K_X))
= H^0 (S^2(L) ( - K_X)) \subset H^0 (S^2(\Omega^1_X) ( - K_X)) .$$

Observe by the way that, if $ L \neq L'$, where we
set $L' : = f^* \Omega^1_B (
\SSS)$, we get a non trivial homomorphism $ L' \ra {\mathcal{ I}}_Z L
(-\Delta)$, hence  $ L - \De \geq L'$.

Since $ 2 L \equiv K_X + \Delta$,
it follows that,
if $F$ is a general fibre, then (use $K_X \cdot F = 0 =L' \cdot F$)
$$ L \cdot F = \Delta \cdot F = 0,$$
hence the effective divisor 	$\Delta$ is contained in a finite
union of fibres.

The first candidate we try with is then the choice of $ L = L'= f^*
\Omega^1_B (
\SSS)$.

To this purpose we recall Kodaira's canonical bundle formula:
$$ K_X \equiv  \SSS_m + f^* (\delta) =  \sum_{i=1}^r
(n_i -1 ) F'_i + f^* (\delta), \ deg (\delta) = \chi(X) -2 + 2b,$$
where $b$ is the genus of the base curve $B$.

Then $H^0 (\hol_X (2 L' - K_X))= H^0 (\hol_X (f^* (2 K_B - \delta)
+ 2 \SSS  - \SSS_m)$, and we search for an effective divisor linearly
equivalent to
$$f^* (2 K_B - \delta)
+ 2 \SSS   - \SSS_m =f^* (2 K_B - \delta)
+ 2 \hat{\SSS}   + \SSS_m .$$

We claim that  $H^0 (\hol_X (2 L' - K_X))= H^0 (\hol_X (f^* (2 K_B -
\delta) )$:
it will then  suffice to have examples where  $|2 K_B - \delta | \neq
\emptyset.$

{\em Proof of the claim}

It suffices to show that
$f_* \hol_X(2
\hat{\SSS}   + \SSS_m) = \hol_B $.
Since the divisor $2 \hat{\SSS}
+ \SSS_m$ is supported on the
singular fibres,
and it is effective,
we have to show that,
    for each singular fibre $F_p =
\sum_{i=1}^{h_p} m_i  C_i$,
neither $2 \hat{\SSS}_p \geq F_p$
nor
${\SSS_m}_{,p} \geq F_p$.

The latter case is obvious since
${\SSS_m}_{,p} = (n_p - 1) F'_p <
F_p = n_p  F'_p $.

In the former
case, $2 \hat{\SSS}_p = \sum_{i=1}^{h_p} 2 (m_i-1)  C_i$,
but it is
not possible that $\forall i$ one has $  2 (m_i-1) \geq m_i $,
since
there is always an irreducible curve $C_i$ with multiplicity $ m_i =
1$.

\QED for the claim

Assume that the elliptic fibration is not
a product (in this case
there is no special
tensor  with vanishing
determinant): then the irregularity of $X$ equals
the genus of $B$,
whence our divisor on the curve $B$ has
degree equal to $ 2b-2 - (1-b
+ p_g(X)) = 3 b -3 - p_g. $

Since $\chi (X) \geq 1$,  $p_g : = p_g
(X) \geq b$, and there exist
an elliptic surface $X$ with any $ p_g
\geq b$ (\cite{qed}).

Since any divisor on $B$ of degree $\geq b$ is effective, it suffices to choose
$  b \leq p_g \leq 2b-3$ and we get a special tensor with trivial determinant,
provided that $b \geq 3$.

Take now a Jacobian elliptic surface in Weierstrass normal form
$$ Z Y^2 - 4 X^3 - g_2 XZ^2 - g_3 Z^3 = 0,$$
where $ g_2 \in H^0 (\hol_B (  4 M))$, $ g_3 \in H^0 (\hol_B (  6 M))$,
and assume that all the fibres are irreducible.

Then the space of special tensors corresponding to our choice of $L$
corresponds to the vector space $H^0 (\hol_B ( 2 K_B - \delta)) =
H^0 (\hol_B (  K_B - 6 M))$. It suffices now to take a hyperelliptic
curve $B$ of genus $ b = 6 h + 1$, and, denoting by $H$ the hyperelliptic
divisor, set $ M : = h H$, so that $K_B - 6 M \equiv 0$
and we have $h^0 (\hol_X (2L-K_X)) =1$. We leave aside for the time
being the question whether the surface $X$ admits a unique special tensor.

\qed

Already in the introduction, we posed the following

{\bf Question. \ }{\em  Let $X$ be a surface with $ q(X)= 0$
and satisfying {\em (1*)}  and {\em (3*)}: is then $X$ strongly
uniformized by the bidisk?}

Concerning the above question, recall the following

      \begin{df}
$  \Gamma \subset \operatorname{Aut} (\HH^n) $ is said to be reducible
      if there exists a subgroup of finite index $\Ga^0 < \Ga$
    such that
$\gamma (z_1, ..., z_n) =( \gamma_1(z_1),..., \gamma_n(z_n))$ for every
$\gamma \in \Gamma^0$) and a decomposition
$\HH^n = \HH^k \times \HH^h$ (with $ h >0$)
     such that the action of $\Gamma^0$  on  $\HH^k $ is
properly discontinuous.
\end{df}

For $n=2$ there are  only two alternatives:

    \begin{oss}

Let  $\Gamma \subset \operatorname{Aut} (\HH^2) $ be a discrete
cocompact
subgroup acting freely
    and  let $X = \HH^2/ \Gamma $.
Then

    \begin{itemize}

    \item $\Gamma$ is reducible if and only
if $X$ is isogenous to a product
of curves, i.e., there is a finite
group $G$ and two curves of genera at least 2
such that
     $X \cong C_1 \times C_2 / G$. Both cases   $q(X) \neq 0$,
$q(X) = 0$ can occur here.

     \item $\Gamma$ is irreducible: then $q(X) = 0$
( this result holds in all dimensions and is a well-known result of Matsushima
     \cite{Ma}).
     \end{itemize}
     \end{oss}

\bigskip

\section{Other surfaces whose universal cover is a product of curves}

For the sake of completeness, using the Enriques classification of surfaces,
we give here a characterization of the K\"ahler surfaces $S$ whose universal
cover is a product of curves, other than $\PP^1 \times \PP^1$ or $ \HH
\times \HH$, which was treated in  section 3. We already mentioned in
the introduction the following theorem.

{\bf  Theorem \ref{surfaces}}
{\em Let $S$ be a smooth compact K\"ahler surface $S$.
Then the universal cover of $S$ is biholomorphic to
\begin{enumerate}
\item
$\PP^1 \times \C$ $\Leftrightarrow$ $ P_{12} = 0$, $q=1$,
$ K_S^2 = 0$.
\item
$\PP^1 \times \HH$ $\Leftrightarrow$ $P_{12} = 0$,
$q=g \geq 2$,
$ K_S^2 = 8 (1-q).$
\item
$\C^2$ $\Leftrightarrow$ $P_{12} = 1$, $ q= 1$ or $q=2$, $K_S^2 = 0$.
\item
$\C \times \HH$ $\Leftrightarrow$ $P_{12} \geq 2$, $ e(S) = 0.$

\end{enumerate}}

\begin{proof}
We consider the several possible cases separately:

1) $\PP^1 \times \C$: by proposition \ref{fibrebundle} these
are the $\PP^1$-
bundles over an elliptic curve. They are  characterized for instance by the
properties
$P_{12} = 0$, which implies that the surface is ruled,
$q=1$, which implies that it is ruled over an elliptic curve,
and $ K^2 = 0$, which implies that all the fibres are smooth, hence we have
a $\PP^1$ bundle.

2) $\PP^1 \times \HH$: these are the $\PP^1$-bundles over a curve $B$ of
genus
$g \geq 2$, hence characterized for instance by the properties $P_{12} = 0$,
$q=g \geq 2$,
$ K^2 = 8 (1-q).$ The argument is here identical to the one given above.

3) $\C^2$: these, by the celebrated theorem of Enriques-Severi
and Bagnera- de Franchis, are the tori or
the hyperelliptic surfaces, characterized
(see for instance  \cite{catcime} page 65), by the properties: $P_{12}
= 1$, $ q= 1$ or
$q=2$,
$K^2 = 0$ (more precisely, $p_g = 1$,  $q=2$, $K^2 = 0$ for tori,
$P_{12} = 1$, $ q= 1$, $K^2 = 0$ for the hyperelliptic surfaces).

4) $\C \times \HH$: in this case, by the same argument as in
proposition \ref{fibrebundle}, the action of $\ga \in \Ga$
is as follows:
$$ (z, \tau) \mapsto ( a_{\ga}(\tau) z +  b_{\ga}(\tau), f_{\ga}(\tau) ),$$
since for fixed $\tau$ we get an automorphism of $\C$.

The cocycle $a_{\ga}(\tau)$ induces a line bundle $L$ which is trivial
on the  leaves $ F_{ \tau} : = (\C \times \{ \tau\}) / \Ga$,
and its dual yields a subbundle of the tangent bundle of $S$.

Moreover, the canonical divisor $K_S$ corresponds to
the cocycle $a_{\ga}(\tau) \cdot \frac {\partial }{\partial \tau}
f_{\ga}(\tau) .$   Therefore the canonical divisor is also trivial on
the leaves
$ F_{ \tau}$, and the extension class of
$$ 0 \ra   \hol_S (K_S - L) \ra \Omega^1_S \ra L \ra 0$$
is given by a group cocycle involving only the function $\tau$.

If the action of $\Ga$ on $\HH$ is properly discontinuous,
then $\HH / \Ga$ is a compact complex curve $B$, and the fibres
of $ f : S \ra B$ are elliptic curves.
There exists  an \'etale cover $S'$ of $S$, such that $S'$ admits an elliptic
fibration with smooth fibres onto a compact complex curve $B'$
of genus at least $2$, hence this is an elliptic bundle (the period map is
constant and $S'$ is K\"ahler).

If the action is not properly discontinuous, then the leaves $ F_{ \tau}$
are not compact. The sections of multiples of the canonical divisor yield
bounded functions on the leaves, hence by Liouville's theorem
these are constant. Since the leaves are  not compact,
the conclusion is that the Kodaira dimension of $S$ is negative or
zero. It cannot be negative, else the universal cover would contain
a family of $\PP^1$'s. If the Kodaira dimension is zero,
we know by surface classification that either the universal cover is $\C^2$
or the fundamental group has order at most two, and in all cases we have
derived a contradiction.

Hence we concluded that our surfaces $S$ are the
elliptic quasi- bundles
$S$ over a curve $B$ of genus $g \geq 2$; more precisely,
these are the
quotients of a product $(E \times C)/G$, where $E$ is an elliptic curve,
$C$ is a
curve of genus
$g' \geq 2$, and $G$ is a finite
group acting diagonally on the product
$E\times C$. These are characterized then by the properties:
$P_{12} \geq 2$, $ e(S) = 0$.

In fact $P_{12} \geq 2$ ensures that
the Kodaira dimension is $ \geq 1$,  a surface of
general type has $e(S) \geq 1$ , whereas for an elliptic fibration $e(S) = 0$
holds if and only if we have a quasi-bundle, i.e., all the fibres are either
smooth or multiple of a smooth curve.

Since $K^2_S = e(S) = 0$, then $\chi (S) = 0$, and Kodaira's canonical bundle
formula says that $K_S$ is the pull back of a $\Q$-divisor on the
base curve $B$ of degree equal to the degree of $K_B +  \sum_{i=1}^r
(n_i -1 ) F'_i $. This means that the base orbifold
is of hyperbolic type, and by the fundamental exact sequence
$ \pi_1(E) \ra \pi_1(S) \ra \pi_1^{orb}(B) \ra 0$ (see \cite{cko} and also
chapter 5 of \cite{catcime}), the universal cover of $S$ is the product
$\C \times \HH$.

    \end{proof}

\section{3-dimensional K\"ahler manifolds  whose
universal cover is  the  polydisk}

In this section we are going to prove theorem \ref{3dim}.

Let $X$  be a smooth compact   K\"ahler manifold of general type of
dimension $3$.
     Assume that the canonical divisor $K_X$ is ample and consider
    the canonical K\"ahler-Einstein metric provided by the theorem of Aubin and
Yau (cf. \cite{Yau}).

    As shown in the introduction, if  $X$ admits a special tensor
    $\omega
\in H^0(X, S^3\Omega^1_X (-K_X))$,  then
by \cite[p.272]{yau1} and  \cite[p.479]{yau}   (see also \cite[p.10]{vz})
$\omega$ induces on the  tangent bundle $T_X$ a homogeneous
    hypersurface  $F_X$ of relative
degree $3$ which is parallel with respect to  the
Levi-Civita connection associated to the K\"ahler-Einstein metric.

    In particular,  taking a point $x\in X$,
    and considering the projectivized tangent bundle,
    we obtain
    a cubic curve $C_x \subset \PP({T_X,x}) \cong \PP^{2}$,
invariant for the action
    of the  holonomy.

    By the  theorem of  De Rham,  the universal cover  $ \tilde{X}$
    splits as a product of irreducible factors,
$\tilde{X} = \tilde{X}_1 \times \tilde{X}_2 \times \cdots \times \tilde{X}_k $
with
$\operatorname{dim}(\tilde{X}_i) = n_i$.
The restricted holonomy group also splits as
$H = H_{1} \times H_{2} \times \cdots \times H_k $, where the action
of $H_i$ on $T_{\tilde{X}_{i}, x_{i}}$ is irreducible
($x_i \in\tilde{X}_i$ being an arbitrary point).

Moreover by the classical theorem of
    Berger-Simons either $H_i \cong U(n_i)$ or  $H_i$ is
    the holonomy of an irreducible Hermitian
symmetric space of rank $> 1$.

    The idea of our proof consists in pointing out    how the 
existence of such a
cubic projective curve
    (possibly singular or reducible)
forces a complete splitting for the action of the holonomy group
    (i.e., it implies the isomorphism   $H \cong U(1)^3$).
Consequently we obtain that $\tilde{X} \cong(\HH)^3$.

\begin{proof}  of theorem \ref{3dim}.

Let $X$  be a smooth  K\"ahler manifold of general type of dimension
$3$, with $K_X$ ample.
    Fix a point $x\in X$ and let $\omega
\in H^0(X, S^3\Omega^1_X (-K_X))$  be a non zero section. Then
$\omega$ induces a projective cubic  curve
$C_x \subset \PP(T_{X,x}) \cong \PP^{2}$  invariant for the action of
the (restricted) holonomy $H$.

In particular $C_x$ is invariant for the action of the  minimal linear
algebraic group which contains $H$, and which we denote
by $\hat{H}$ .  Observe that $\hat{H}$ is connected.

On the other side,  by the  description given above,   we have
$H = H_{1} \times H_{2} \times \cdots \times H_k $,  where    either
$H_i \cong U(n_i)$ or  $H_i$ is
    the holonomy of an irreducible Hermitian symmetric space of rank $> 1$.

Let $\operatorname{Lin} (C_x)$   be the linear algebraic group
of projectivities leaving  $C_x$ invariant.
We shall analyse all the possible cases for $C_x$, including the study of
its  singularities and the description of $\operatorname{Lin}
(C_x)$,
    keeping in mind that  we have $\PP (\hat{H} )\subset \operatorname{Lin}
(C_x)$.

\subsubsection*{(a)  $C_x$ irreducible and smooth. }
In this case $\operatorname{Lin} (C_x)$ is finite, which contradicts $\PP
(\hat{H} )
\subset \operatorname{Lin} (C_x)$, since  dim $\hat{H}$ is at least $3$.

\subsubsection*{(b)  $C_x$ irreducible with a node $p$. }
    In this case   $\hat{H}$ fixes the node  $p$  and  the pair of
tangent lines of
$C_x$ at $p$.  Since  $\hat{H}$  is connected,  it fixes both tangent
lines.

     Therefore $H$ fixes the  point $p$ and a line  $L$ through $p$, i.e.
$ H$ fixes a flag.
    Since $H$ is a subgroup of the unitary subgroup   it acts diagonally
for a suitable unitary basis, hence we conclude that $ H = U(1)^3$.

Therefore there exists an \'etale covering $X'$ of $X$ such that
    $T_{X'}$  decomposes as the direct sum of   3   line bundles
(the eigenbundles of the action), and  the universal
cover of $X$ is biholomorphic to $\HH^3$.

\subsubsection*{(c)  $C_x$ irreducible with a cusp $p$. }  In this case
we can choose coordinates on $\PP^2$ so that $p=(1:0:0)$ and  on the
affine chart
$ \{x_0=1\}$  the curve $C_x$ is parametrized by $t \mapsto (1,t^2, t^3)$.

    Now we have :
$\C^{\ast} \cong \operatorname{Lin} (C_x)  $ and in the
affine chart $ \{x_0=1\}$
$\lambda \in \C^{\ast}$ yields the automorphism
$$ \begin{array}{rll}  C_x  & \rightarrow & C_x \\
(1,t^2,t^3) & \mapsto &  (1, \lambda^2t^2, \lambda^3t^3)
\end{array}
$$
Whence even in this case  the action of $\hat{H}$
is diagonal  and we  conclude as before.

    \subsubsection*{(d)  $C_x$ decomposes as the union of a line $L$
and an irreducible conic $Q$. }
    In this case $\hat{H}$ fixes the intersection set $L \cap Q $,
which consists of one or two points. By connectedness of  $\hat{H}$,
$\hat{H}$  fixes a point $ P \in L$ and
the line $L$, and we conclude as before.

    \subsubsection*{(e)  $C_x$ decomposes as the union of a double
line $2L_1$ and a line $L_2$. }
    In this case $\hat{H}$ fixes the point $L_1 \cap L_2 = \{p\}$ and
the line $L_2$ and we  are done.

     \subsubsection*{(f)  $C_x$ decomposes as the union of 3 distinct
lines $C_x= L_1\cup L_2 \cup L_3$. }
     There are two possibilities: the three lines are concurrent
in the same point $p$ or
     $ L_1\cap L_2 \cap L_3 =\emptyset$
and there are three singular points $p_{ij}=L_i\cap L_j$ ($ 1\leq i<j\leq3$).

    In both   cases, since $\hat{H}$ is connected it fixes each
singular point and
each line. Hence there is a  flag fixed by $\hat{H}$ and we are done.

     \subsubsection*{(g)  $C_x$ decomposes as   a triple line $3L$. }
    We are going to show that this case cannot happen.

     Assume the contrary and
consider  the line subbundle ${\mathcal L} \subset  \Omega_X^1
$ corresponding to $L$. We have a section
$$  \hol_X \ra \hol_X(3{\mathcal L} -K_X)  \subset S^3 \Omega^1_X (-K_X)
$$ (indeed, cf. \cite{yau} or \cite{vz}, this section has no zeros).

Therefore we  have $3 {\mathcal L}  \equiv K_X  + D$, with $D$ effective
( in fact $D$ is a trivial divisor). This
in particular    implies  $\mathcal L$  big because $K_X$  is ample by
our assumption.
     This assertion, as in the proof of theorem  \ref{unibidisk},
contradicts
     the theorem of  Bogomolov (cf. \cite{Bog}).

     \hfill\break
     Conversely, if
    $X \cong  \HH  \times \HH \times \HH  / \Gamma$, with  $\Gamma$ a
cocompact discrete  subgroup of $\operatorname{Aut} (\HH \times \HH \times
\HH) $ acting freely, then by \cite {siegel}  it is immediately
seen that $K_X $ is
ample  and by Prop. \ref{nec}
$X$ admits a semi special tensor.

     \end{proof}

\section{4-dimensional K\"ahler manifolds of general type
with a special tensor
whose universal cover is not  a product of curves}

One of the consequences of the theorem of Berger-Simons is that an
irreducible K\"ahler   manifold $X$  of dimension $n$ and with $K_X$ ample
(irreducible in the sense of  De Rham's  theorem) has  as holonomy group
a proper subgroup
$H \subset U(n)$  if and only if $\tilde{X}$ is a Hermitian symmetric space
     of rank $\geq 2$ (see \cite{yau1}, and especially \cite{vz}, 1.4 and 1.5).

Since we are interested in the case where  $K_X$ is ample
we look for the Cartan realization of a
Hermitian symmetric space of noncompact type as a
bounded complex symmetric domain.

We shall find first such a bounded symmetric domain such that
it has a holonomy invariant hypersurface of degree $n$,
and then we shall apply
     the classical result of Borel  on complex analytic Clifford-Klein
     forms.
A complex analytic Clifford-Klein
     form is simply a compact quotient $ X = \tilde{X} / \Ga$,
where the group $\Ga$ acts freely (thus $X$ is a projective manifold
with ample canonical bundle).

Borel's theorem (cf. \cite{bo63}) states that any
    bounded symmetric domain $\tilde{X}$ of dimension $n$  admits  infinitely
many
    compact complex analytic Clifford-Klein
     forms,
     whose arithmetic genus $ 1 - \chi(X)$ can be  arbitrarily large in absolute
value.

     We shall prove Theorem \ref{4dim}  considering a Clifford-Klein
     form  $X$ associated to the  noncompact  Hermitian symmetric space
    of complex dimension 4
    $\tilde{X} : =SU(2,2) / S (U(2) \times U(2))$. In higher dimensions, it
clearly suffices to take the product of such a projective manifold $X$ with
$n-4$ projective curves $ C_1, \dots C_{n-4}$ of genus at least $2$.

\hfill\break
\Proof of Theorem \ref{4dim}

Let $\tilde{X}= SU(2,2) / S (U(2) \times U(2))$.
$\tilde{X}$ is a noncompact  Hermitian symmetric space
     of dimension 4 and rank 2.
Recall that a  $4 \times 4$  matrix $g \in S U(2,2)$ can be written as
     $$g= \begin{pmatrix}
A & B \\
C & D
\end{pmatrix}$$   where   $\det (g)=1$ and $A,B, C,D$
are  $2\times 2$ complex matrices
satisfying
$$ (\star) \ \
    ^{t}  \overline{A} \cdot A - ^{t}  \overline{C}\cdot
C=\operatorname{Id}  ;  \
     ^{t}  \overline{B} \cdot B - ^{t}  \overline{D}\cdot D= -
\operatorname{Id}  ;  \
    ^{t}  \overline{B} \cdot A - ^{t}  \overline{D}\cdot C =0 , \
    $$
    whereas the subgroup $S (U(2) \times U(2))$ can be identified with
the matrices  of the form
    $$\begin{pmatrix}
A & 0 \\
0 & D
\end{pmatrix} \ \ \  (\mbox{with } \  A, D \in U(2) \ , \ \det (A)
\cdot \det (D) =1).  $$

    Let $\mathfrak{su} (2,2)$ be the Lie algebra   of $SU(2,2)$.
The Cartan decomposition
    $\mathfrak{su} (2,2) =  \mathfrak{k} \oplus \mathfrak{p}$
can be written down  explicitly by means   of
$$ \mathfrak{p}  \cong
    \begin{pmatrix}
0 & B \\
^{t}  \overline{B} & 0
\end{pmatrix} \  ,  \  \mathfrak{k}  \cong
    \begin{pmatrix}
A & 0 \\
    0 & D
\end{pmatrix} \ \  (\mbox{with } \  ^{t}  \overline{A} = -A  \ ,  \
^{t}  \overline{D} = -D ) $$
and  for $x\in \tilde{X}$ we have a canonical
isomorphism   $ \mathfrak{p}  \cong T_{X,
x}$.

The  holonomy action
    coincides with the adjoint representation of  $S (U(2) \times U(2))$
on $\mathfrak{p}$, given
    for every matrix   $M= \begin{pmatrix}
A & 0 \\
0 & D
\end{pmatrix}  \in S (U(2) \times U(2))$   by  the map
$\operatorname{Ad}_{M} :    \mathfrak{p}  \rightarrow   \mathfrak{p}  $
described by
$$ \begin{array} {rccc}
    & \begin{pmatrix}
0 & B \\
^{t}  \overline{B} & 0
\end{pmatrix} & \mapsto
\begin{pmatrix}
0 & A\cdot B \cdot D^{-1} \\
    D \cdot ( ^{t}  \overline{B})  \cdot  A^{-1}  & 0
\end{pmatrix}
\end{array}
$$

    Let us now consider  the Cartan realization of
     $\tilde{X}$.
It    is obtained by the Siegel domain in the space of $2\times 2$ matrices
     $M_{2,2}(\C)$ (see \cite[p.527]{hel}):
     $$ X \cong   \{ Z \in  M_{2,2}(\C) \  :  \  \operatorname{Id} - ^{t} Z
    \cdot \overline{Z} >0\}$$
     and
    the action of $SU(2,2) $  on $X$ is given by:
$$ Z \mapsto (AZ+B)\cdot (CZ+D)^{-1}  $$
Considering  the tangent space at  $0$,  the action of $S (U(2)
\times U(2))$  on an "infinitesimal" $ 2 \times 2$ matrix  $Z$ becomes
$$ Z \mapsto AZ D^{-1}  $$
     and  in particular  we recover  the above description of
     the adjoint representation of  $S (U(2) \times U(2))$.

    Notice that, since $\det (A) \cdot \det (D) =1$,
     we have $\det AZ D^{-1}  = \det(A)^2 \cdot \det Z$.
     This exactly means that   the determinant is a semi-invariant for
the action of
$S (U(2) \times U(2))$  on $T_{X,0}$.

Therefore,  identifying $T_{X,0}$  with $M_{2,2}$, and   considering
     the   projectivized  tangent bundle   at $0$,   $\PP(T_{X,0})\cong \PP^3$,
    $\{   \det(Z) =0\}$  defines a quadric surface,
invariant for the action of $S (U(2) \times U(2))$, and of course
     we obtain an invariant  quartic  projective surface
    given by $\{ Z \in M_{2,2} \
:   (\det(Z))^2 =0\}$.

    Applying now the theorem of Borel cited above
we obtain a compact complex analytic Clifford-Klein
     form $ X \cong \tilde{X} / \Gamma$ of $\tilde{X}$.
We shall exhibit a semispecial tensor   $ \tilde{\omega} $ on $\tilde{X}$
which will descend to $X$ yielding a semispecial tensor. Since $\tilde{X}$
is irreducible, our proof will be complete.

We want to show how this invariant surface defines a special tensor.

Write, for $ \ga \in \Ga$,
$$ \ga (Z)  =  (AZ+B)\cdot (CZ+D)^{-1}  \Leftrightarrow
\ga (Z) \cdot (CZ+D) = (AZ+B).$$
Differentiating the above equality, we obtain
$$ d \ga (Z) \cdot (CZ+D) = ( A - (AZ+B)\cdot (CZ+D)^{-1} C )\cdot d Z.$$

Taking determinants,  we obtain
$$ det (d \ga (Z)) \cdot det (CZ+D) = det ( A - (AZ+B)\cdot (CZ+D)^{-1} C )
\cdot
det (d Z)=$$
$$ = det ((CZ+D)^{-1}) det (C) det (AC^{-1} D - B) \cdot
det (d Z). $$
Observe now that, setting $^*B : = ^t \bar{B}$, equations $(\star )$ yield
   $$det
(AC^{-1} D - B) = det (^*{B^{-1}} ^* D D - B) = det (^*{ B^{-1}}).$$

An easy calculation using the above equations yields then
$det (C) det (AC^{-1} D - B) = det (A) \det (^* D)^{-1} = det (A) \det ( D)$.

If we restrict to the isotropy subgroup $H = S (U(2) \times U(2))$,
    we get  $det
(A) \cdot det ( D) = 1$. We have now a character of the group which is trivial
on $H$. This character is then trivial since the homogeneous domain is
contractible, whence the group  $G : = S (U(2,2))$ is homotopically equivalent
to $H$.

Since finally $det
((CZ+D)^{-4})$ is the Jacobian determinant of the transformation
$\ga$,
$ \tilde{\omega} : = det  (d Z)^2$ is a $\Ga$-invariant section of
$H^0(\tilde{X}, S^n\Omega^1_{\tilde{X}} (-K_{\tilde{X}}))$,
thus a special tensor which descends to $X$.

\qed

\hfill\break  {\bf Acknowledgements.}

These research   was performed in the realm  of the D.F.G. Programs:
       SCHWERPUNKT "Globale Methoden in der komplexen Geometrie",
and
FORSCHERGRUPPE 790
`Classification of algebraic surfaces and compact complex
manifolds'.

The second author thanks the Universit\"at Bayreuth for its warm
hospitality
     in the months of november and december 2006 (where the research was
begun) and the DFG  for supporting his  visit. The first author
gratefully acknowledges the   hospitality of  the Centro De Giorgi,
where the final part  of the present article was written.

We would like to thank Eckart Viehweg and Kang Zuo for pointing out
their clear explanation (\cite{vz}) of some aspects of
Yau's uniformization theorem in terms of stability of symmetric powers of
the cotangent bundle, J\"org Winkelmann for pointing out the result of
Fornaess and Stout,  Antonio J. Di Scala for answering some queries
of the second author, Fedor Bogomolov for an interesting conversation.

\noindent {\bf Authors' addresses:}

\bigskip

\noindent Prof. Fabrizio Catanese\\ Lehrstuhl Mathematik VIII,
Universit\"at Bayreuth, NWII\\
       D-95440 Bayreuth, Germany \\  e-mail: Fabrizio.Catanese@uni-bayreuth.de

\bigskip

\noindent Marco Franciosi \\ Dipartimento di Matematica Applicata "U.
Dini", Universit\`a di Pisa \\
      via Buonarroti 1C, I-56127, Pisa, Italy \\  e-mail:
franciosi@dma.unipi.it

\end{document}